# Asymptotic behavior of summation functions from bounded arithmetic multiplicative functions

Victor Volfson

ABSTRACT  The paper considers estimates for the asymptotics of summation functions of bounded multiplicative arithmetic functions. Several assertions on this subject are proved and examples are considered.

# 1. INTRODUCTION

In the general case, an arithmetic function is a function defined on the set of natural numbers and taking values on the set of complex numbers. The name arithmetic function is due to the fact that this function expresses some arithmetic property of the natural series.

A multiplicative arithmetic function is an arithmetic function for which:

$$g(m_1 m_2) = g(m_1) g(m_2),$$

where $m_1$, $m_2$ are coprime numbers.

Let a natural number $m$ have a canonical decomposition $m = p_1^{a_1} ... p_t^{a_t}$, where $p_i$ is a prime number and $\alpha_i$ is a natural number. Then, it is holds for the multiplicative arithmetic function:

$$g(m) = g(p_1^{a_1} ... p_t^{a_t}) = g(p_1^{a_1}) ... g(p_t^{a_t}) = \prod_{p^\alpha | m} g(p^\alpha). \qquad (1.1)$$

A strongly multiplicative arithmetic function is an arithmetic function for which is true - $g^*(p^\alpha) = g^*(p)$. Therefore, based on (1.1), it is executed for a strongly multiplicative arithmetic function:

$$g^*(m) = g^*(p_1^{a_1} ... p_t^{a_t}) = \prod_{p | m} g^*(p). \qquad (1.2)$$

The summation function $S$ from an arithmetic function $f(m), m = 1,...,n$ is called:

$$S(n) = \sum_{m \leq n} f(m). \qquad (1.3)$$

Let's look at the summation function from the restricted arithmetic function, i.e. $|f(m)| \leq C$ at $m = 1,...,n$.

If $C > 1$, then take $g(m) = f(m) / C$, for $m = 1,...,n$, for which $|g(m)| \leq 1$. Therefore, if there is $\lim_{n \to \infty} \dfrac{\sum_{m \leq n} f(m)}{n}$, then it is equal to:

$$\lim_{n\to\infty}\frac{\sum_{m\leq n}f(m)}{n}=\lim_{n\to\infty}\frac{\sum_{m\leq n}Cg(m)}{n}=C\lim_{n\to\infty}\frac{\sum_{m\leq n}g(m)}{n}. \qquad (1.4)$$

Further, if $g(m), m=1,...,n$ is a multiplicative arithmetic function and $|g(m)|\leq 1$, then according to the Wirzing theorem [1]:

$$\lim_{n\to\infty}\frac{\sum_{m\leq n}g(m)}{n}=\prod_p(1-\frac{1}{p})(1+\frac{\sum_{\alpha\geq 1}g(p^\alpha)}{p^\alpha}). \qquad (1.5)$$

Based on (1.4) and (1.5) we get:

$$\lim_{n\to\infty}\frac{\sum_{m\leq n}f(m)}{n}=C\lim_{n\to\infty}\frac{\sum_{m\leq n}g(m)}{n}=C\prod_p(1-\frac{1}{p})(1+\frac{\sum_{\alpha\geq 1}g(p^\alpha)}{p^\alpha}). \qquad (1.6)$$

Taking into account (1.6), we will further consider only the case with $|g(m)|\leq 1, m=1,...,n$.

An example of a summation function from a bounded multiplicative arithmetic function is the Mertens function:

$$M(n)=\sum_{m\leq n}\mu(m), \qquad (1.7)$$

where $\mu(m)$ is the Möbius function, where the Möbius function $\mu(m)=1$ if the natural number m has an even number of prime divisors of the first degree, $\mu(m)=-1$ if the natural number m has an odd number of prime divisors of the first degree, and $\mu(m)=0$ if the natural number m has prime divisors not only of the first degree.

Finding the asymptotic behavior of summation functions from bounded multiplicative arithmetic functions is an important problem in number theory [1-3]. This work is devoted to solving this problem.

2. ASYMPTOTIC BEHAVIOR OF SUMMATION FUNCTIONS FROM A BOUNDED MULTIPLICATIVE ARITHMETIC FUNCTIONS

The assertion was proved in [4].

Assertion 1

Let $|g(m)| \leq 1, m = 1,...,n$ is a real multiplicative arithmetic function, then the relation is true for the summation function:

$$S(n) = \sum_{m \leq n} g(m) = n \prod_p (1 - \frac{1}{p}) \sum_{\alpha \geq 0} \frac{g(p^\alpha)}{p^\alpha} + o(n),$$

where the infinite product is considered to be 0 if the product diverges.

Let us consider several consequences of Assertion 1.

Assertion 2

Let $|g(m)| \leq 1, m = 1,...,n$ is a real multiplicative arithmetic function, then if the series $\sum_p \frac{1 - g(p)}{p} < \infty$ (converges), then the infinite product converges and the asymptotic of the summation function is equal to:

$$\sum_{m \leq n} g(m) = n \prod_p (1 - \frac{1}{p}) \sum_{\alpha \geq 0} \frac{g(p^\alpha)}{p^\alpha} + o(n). \tag{2.1}$$

If the series $\sum_p \frac{1 - g(p)}{p} = \infty$ (diverges), then the infinite product (2.1) tends to zero, and the asymptotic of the summation function is equal to:

$$\sum_{m \leq n} g(m) = o(n). \tag{2.2}$$

Proof

We write the infinite product in the form:

$$\prod_p (1 - \frac{1}{p}) \sum_{\alpha \geq 0} \frac{g(p^\alpha)}{p^\alpha} = \prod_p (1 - \frac{1 - g(p)}{p} - \frac{g(p) - g(p^2)}{p^2} - ...) \tag{2.3}$$

Taking into account that $|g(p^\alpha)| \leq 1, \alpha \geq 1$, then each term of the infinite product (2.3) tends to 1 as at $p \to \infty$. Therefore, starting from some $N$, each term of the infinite product is positive and the product can be taken as a logarithm:

$$\ln(\prod_{p\geq N}(1-\frac{1-g(p)}{p}-\frac{g(p)-g(p^2)}{p^2}-...)) = \sum_{p\geq N}\ln(1-\frac{1-g(p)}{p}+O(\frac{1}{p^2})) = -\sum_{p\geq N}\frac{1-g(p)}{p}+O(\sum_{p\geq N}\frac{1}{p^2}) .(2.4)$$

Since the series $\sum_p \frac{1}{p^2}$ converges, if the series $\sum_p \frac{1-g(p)}{p}$ converges, then the series $\sum_{p\geq N}\frac{1-g(p)}{p}$ converges and the series (2.4) converges. Therefore, in this case, the infinite product converges:

$$\prod_{p\geq N}(1-\frac{1-g(p)}{p}-\frac{g(p)-g(p^2)}{p^2}-...) \qquad (2.5)$$

and also the infinite product (2.3). Therefore, by Assertion 1, (2.1) holds.

If the series $\sum_p \frac{1-g(p)}{p} = \infty$ (diverges), then the series (2.4) diverges to $-\infty$, and the infinite products (2.5) and (2.3) tend to zero. Therefore, by Assertion 1, (2.2) holds.

An example for assertions 1 and 2.

Let's find the asymptotic of the Mertens function $M(n) = \sum_{m\leq n}\mu(m)$.

It is known that $\mu(m), m = 1,...,n$ is a multiplicative arithmetic function and $|\mu(m)|\leq 1$, therefore, assertion 1 can be used to find the indicated asymptotic.

Since $\mu(p) = -1, \mu(p^\alpha) = 0$ at $\alpha > 1$, then

$$\sum_{m\leq n}\mu(m) = n\prod_p(1-\frac{1}{p})(1-\frac{1}{p}+\frac{0}{p^2}+...)+o(n) = n\prod_p(1-\frac{1}{p})^2 + o(n) .$$

Having in mind that $\prod_{p\leq n}(1-\frac{1}{p})^2 > 0$, therefore, we can take the logarithm of the product:

$$\ln(\prod_{p\leq n}(1-\frac{1}{p})^2) = \sum_{p\leq n}\ln(1-\frac{1}{p})^2 = -2\sum_{p\leq n}\frac{1}{p}+O(\sum_{p\leq n}\frac{1}{p^2}) = -2\ln\ln n + O(1).$$

Therefore, the asymptotic of the product will be equal to:

$$\prod_{p\leq n}(1-\frac{1}{p})^2 \sim \frac{e^{-c}}{(e^{\ln\ln n})^2} = \frac{e^{-c}}{(\ln n)^2} .$$

Then the infinite product will be equal to:

$$\prod_p (1-\frac{1}{p})^2 = \lim_{n\to\infty} \prod_{p\leq n}(1-\frac{1}{p})^2 = \lim_{n\to\infty} \frac{e^{-c}}{(\ln n)^2} = 0.$$

Taking this into account, based on Assertion 1, the asymptotic of the Mertens function will be equal to:

$$M(n) = \sum_{m\leq n} \mu(m) = o(n).$$

Now we use Assertion 2.

The series $\sum_p \frac{1-\mu(p)}{p} = \sum_p \frac{1-(-1)}{p} = \sum_p \frac{2}{p}$ - diverges, therefore $\sum_{m\leq n} \mu(m) = o(n)$.

The second proof is shorter.

The zeta function of primes was considered in [5]. It is defined as the following infinite series converging for Re(s) > 1:

$$P(s) = \sum_p \frac{1}{p^s}, \qquad (2.6)$$

where $p$ is a prime number.

The zeta function of primes is related to the Riemann zeta function $\zeta$ as follows [6]:

$$P(s) = \sum_{n=1}^{\infty} \frac{\mu(n)}{n} \ln \zeta(ns). \qquad (2.7)$$

Assertion 3

Let $0 < g(m) \leq 1, m = 1,...,n$ is a real multiplicative arithmetic function, then if the series

$$\sum_p \frac{1-g(p)}{p} < \infty \quad \text{(converges)} \quad \text{and} -1 < \frac{1-g(p)}{p} + \frac{g(p)-g(p^2)}{p^2} + ... + \frac{g(p^{k-1})-g(p^k)}{p^k} + ... < 1,$$

then the asymptotic of the summation function in terms of the zeta function of primes is:

$$\sum_{m\leq n} g(m) = e^{-c} n + o(n), \qquad (2.8)$$

where $c = \sum_{k=2}^{\infty} b_k P(k)$.

The asymptotic of the summation function through the Riemann zeta function under the same conditions is determined by the formula (2.8), where

$$P(k) = \sum_{n=1}^{\infty} \frac{\mu(n)}{n} \ln \zeta(nk). \tag{2.9}$$

Proof

Having in mind that $0 < g(m) \leq 1, m = 1,...,n$, then the infinite product is:

$$\prod_p (1 - \frac{1}{p}) \sum_{\alpha \geq 0} \frac{g(p^\alpha)}{p^\alpha} > 0. \tag{2.10}$$

Based on (2.10), the infinite product can be taken as a logarithm and, using (2.3), we obtain:

$$\ln(\prod_p (1 - \frac{1-g(p)}{p} - \frac{g(p)-g(p^2)}{p^2} - ... - \frac{g(p^{k-1})-g(p^k)}{p^k} - ...)) = \sum_p \ln(1 - \frac{1-g(p)}{p} - \frac{g(p)-g(p^2)}{p^2} - ... - \frac{g(p^{k-1})-g(p^k)}{p^k} - ...)) \tag{2.11}$$

We expand (2.11) in a Taylor series, having in mind that $-1 < X(p) = \frac{1-g(p)}{p} + \frac{g(p)-g(p^2)}{p^2} + ... + \frac{g(p^{k-1})-g(p^k)}{p^k} + ... < 1$:

$$\sum_p \ln(1 - X(p)) = -\sum_p X(p) + \sum_p \frac{X^2(p)}{2} - ... + (-1)^n \sum_p \frac{X^n(p)}{n} + ..., \tag{2.12}$$

where $X(p) = \frac{1-g(p)}{p} + \frac{g(p)-g(p^2)}{p^2} + ... + \frac{g(p^{k-1})-g(p^k)}{p^k} + ... = \frac{a_2}{p^2} + \frac{a_3}{p^3} + ... + \frac{a_k}{p^k} + ...$.

There is no term $\frac{a_1}{p}$ in (2.12), otherwise the series $\sum_p \frac{1-g(p)}{p}$ would diverge, which contradicts the condition.

Based on (2.6), series (2.12) can be written as:

$$\sum_p \ln(1 - \sum_{i=2}^{\infty} \frac{a_i}{p^i}) = -\sum_p (\frac{a_2}{p^2} + \frac{a_3}{p^3} + ...) + \frac{1}{2} \sum_p (\frac{a_2}{p^2} + \frac{a_3}{p^3} + ...)^2 - ... + \frac{(-1)^n}{n} \sum_p (\frac{a_2}{p^2} + \frac{a_3}{p^3} + ...)^n + ... = -\sum_{k=2}^{\infty} b_k P(k). \tag{2.13}$$

Having in mind (2.13), in this case, the infinite product can be written as:

$$\prod_p (1-\frac{1}{p})\sum_{\alpha \geq 0}\frac{g(p^{\alpha})}{p^{\alpha}} = e^{-\sum_{k=2}^{\infty}b_k P(k)}. \qquad (2.14)$$

Based on (2.14) and Assertion 2, we obtain the asymptotic for the summation function in terms of the zeta function of primes:

$$\sum_{m \leq n} g(m) = n e^{-\sum_{k=2}^{\infty}b_k P(k)} + o(n), \qquad (2.15)$$

which corresponds to (2.8).

Having in mind (2.7) and (2.15) we obtain the asymptotic of the summation function in terms of the Riemann zeta function (2.8), (2.9).

Assertion 4

Let it will be $|g(m)| \leq 1, m = 1,...,n$ strongly multiplicative arithmetic function, then the relation for the summation function is true:

$$\sum_{m \leq n} g^*(m) = n\prod_p (1-\frac{1-g^*(p)}{p}) + o(n), \qquad (2.16)$$

where the infinite product is considered to be 0 if the product diverges.

The proof follows from the definition of a strongly multiplicative arithmetic function $g^*(p^{\alpha}) = g^*(p)$ and Assertion 1.

Let's take an example of assertions 3 and 4.

The arithmetic function $g^*(m) = \frac{\varphi^2(m)}{m^2} \leq 1$ is real strongly multiplicative, $g^*(p) = \frac{(p-1)^2}{p^2} = \frac{p^2-2p+1}{p^2} = 1-\frac{2}{p}+\frac{1}{p^2}$ and $\sum_p \frac{1-g^*(p)}{p} < \infty$. Therefore, based on Assertion 4, we get:

$$\prod_p (1-\frac{1-g^*(p)}{p}) = \prod_p (1-\frac{2}{p^2}+\frac{1}{p^3}). \qquad (2.17)$$

All members of the product tend to 1 at $p \to \infty$, so it is greater than zero and this product can be logarithmized:

$$\ln(\prod_p (1 - \frac{2}{p^2} + \frac{1}{p^3})) = \sum_p \ln(1 - \frac{2}{p^2} + \frac{1}{p^3}). \tag{2.18}$$

Let us expand (2.18) in a Taylor series:

$$\sum_p \ln(1 - \frac{2}{p^2} + \frac{1}{p^3}) = -2\sum_p \frac{1}{p^2} + \sum_p \frac{1}{p^3} - \frac{1}{2}\sum_p (-\frac{2}{p^2} + \frac{1}{p^3})^2 + \ldots + \frac{(-1)^{n+1}}{n}\sum_p (-\frac{2}{p^2} + \frac{1}{p^3})^n + \ldots \tag{2.19}$$

Based on (2.19) we get:

$$\sum_p \ln(1 - \frac{2}{p^2} + \frac{1}{p^3}) = -2\sum_p \frac{1}{p^2} + \sum_p \frac{1}{p^3} - 2\sum_p \frac{1}{p^4} + 2\sum_p \frac{1}{p^5} - \frac{19}{6}\sum_p \frac{1}{p^6} + 4\sum_p \frac{1}{p^7} - 2\sum_p \frac{1}{p^8} + \frac{1}{3}\sum_p \frac{1}{p^9} + \ldots \tag{2.20}$$

The values of the zeta function of prime numbers are in the OEIS tables: A085548, A085541, A085964, A08596, so based on (2.20) we get:

$$\sum_p \ln(1 - \frac{2}{p^2} + \frac{1}{p^3}) \approx -0{,}8122. \tag{2.21}$$

Having in mind (2.21) we get:

$$\prod_p (1 - \frac{2}{p^2} + \frac{1}{p^3}) \approx e^{-0{,}8122} = 0{,}4438. \tag{2.22}$$

Taking into account (2.22), based on Assertion 4, we obtain the asymptotic of the summation function:

$$\sum_{m \leq n} g^*(m) = n \prod_p (1 - \frac{1 - g^*(p)}{p}) + o(n) \approx 0{,}4438n + o(n). \tag{2.23}$$

3. ASYMPTOTIC OF THE NUMBER OF NATURAL NUMBERS SATISFYING CERTAIN REQUIREMENTS

The number of natural numbers that meet certain requirements is a summation function of the form $\sum_{m \leq n} 1_\beta(m)$, where $1_\beta(m)$ is an arithmetic indicator function. $1_\beta(m) = 1$, if $m \in \beta$, where $\beta$ is a set of natural numbers that satisfies certain requirements, otherwise $1_\beta(m) = 0$.

We consider hear only the case when the indicator arithmetic function $1_\beta(m)$ is multiplicative.

The simplest multiplicative arithmetic indicator function is the arithmetic function $\varepsilon(m)$ ($\varepsilon(m) = 1$, if $m = 1$ and $\varepsilon(m) = 0$, if $m > 1$). It is clear that $\sum_{m \leq n} \varepsilon(m) = 1$.

Another multiplicative arithmetic indicator function is the arithmetic function $1(m)$ ($1(m) = 1$ for any value $m$). It's clear that $\sum_{m \leq n} 1(m) = n$.

Another example of an indicator from a multiplicative arithmetic function is square-free natural numbers (this is a set of natural numbers $\beta$ whose decomposition into a product of powers of prime numbers contains only the first powers, i.e. $m \in \beta$, if $m = p_1 p_2 ... p_k$).

The number of natural numbers not greater than $n$, which are free from squares, is usually denoted by $Q(n) = \sum_{\substack{m \leq n \\ m = p_1 p_2 ... p_k}} 1(m)$.

Taking into account that $g(m) = 1_\beta(m)$ is a multiplicative function in this case and $|g(m)| \leq 1$, then as $\sum_p \frac{1 - g(p)}{p} = \sum_p \frac{1-1}{p}$ - converges, then we determine the asymptotic using Assertion 2 and the Euler identity.

$$Q(n) = \sum_{m \leq n} g(m) = n \prod_p (1 - \frac{1}{p})(1 + \frac{1}{p} + \frac{0}{p^2} + ...) + o(n) = n \prod_p (1 - \frac{1}{p^2}) + o(n) = \frac{n}{\zeta(2)} + o(n) . \quad (3.1)$$

Based on (3.1), the asymptotic of the number of natural numbers not exceeding $n$, which are not free from squares, i.e. the expansion of which into a product of powers of primes contains not only the first powers (we denote $Q^-(n)$), is defined as:

$$Q^-(n) = n - Q(n) = (1 - \frac{1}{\zeta(2)})n + o(n) = O(n) . \quad (3.2)$$

Let's look at another example of determining the asymptotic of the number of natural numbers that satisfy certain requirements, when the indicator arithmetic function is multiplicative.

This is an arithmetic function $g(m) = powerful(m)$. The value of $powerful(m) = 1$, if $m = 1$ or in the factorization of a natural number $m = p_1^{\alpha_1} p_2^{\alpha_2} ... p_k^{\alpha_k}, \alpha_1 > 1, \alpha_2 > 1, ..., \alpha_k > 1$, otherwise - $powerful(m) = 0$. It is easy to verify that the arithmetic function $powerful(m)$ is real and multiplicative.

It is clear that the natural numbers for which $powerful(m) = 1$ is less than the numbers not free from squares (3.2). How much less?

Let us determine the asymptotic of the number of natural numbers, where $powerful(m) = 1$, using Assertion 1.

Having in mind that $g(p) = 0, g(p^\alpha) = 1, \alpha > 1$ we get:

$$\prod_{p \leq n}(1-\frac{1}{p})\sum_{\alpha \geq 0}\frac{g(p^\alpha)}{p^\alpha} = \prod_{p \leq n}(1-\frac{1}{p})(1+\frac{1}{p^2}+...+\frac{1}{p^\alpha}) = \prod_{p \leq n}(1-\frac{1}{p}+\frac{1}{p^2}) > 0. \qquad (3.3)$$

We take the logarithm of (3.3) and obtain:

$$\ln(\prod_{p \leq n}(1-\frac{1}{p}+\frac{1}{p^2})) = \sum_{p \leq n}\ln(1-\frac{1}{p}+\frac{1}{p^2}) = -\sum_{p \leq n}\frac{1}{p}+O(\sum_{p \leq n}\frac{1}{p^2}) = -\ln\ln n + O(1). \qquad (3.4)$$

Based on (3.4), the asymptotic of the product is:

$$\prod_{p \leq n}(1-\frac{1}{p}+\frac{1}{p^2}) \sim \frac{e^{-c}}{e^{\ln\ln n}} = \frac{e^{-c}}{\ln n}. \qquad (3.5)$$

Taking into account (3.5) the infinite product is:

$$\prod_{p}(1-\frac{1}{p}+\frac{1}{p^2}) = \lim_{n \to \infty}\prod_{p \leq n}(1-\frac{1}{p}+\frac{1}{p^2}) = \lim_{n \to \infty}\frac{e^{-c}}{\ln n} = 0. \qquad (3.6)$$

Based on (3.6) and Assertion 1:

$$\sum_{m \leq n} powerful(m) = o(n). \qquad (3.7)$$

Now let's do this using Assertion 2.

Using that $powerful(p) = 0$, the series $\sum_{p}\frac{1-powerful(p)}{p} = \sum_{p}\frac{1}{p}$ - diverges, therefore:

$$\sum_{m \leq n} powerful(m) = o(n).$$

Naturally, estimate (3.7) is not exact. An accurate estimate can be obtained using complex analysis methods [7]. However, this estimate is suitable for qualitative analysis. Let us compare, for example, the asymptotic estimates for two arithmetic functions (3.2) and (3.7).

4. CONCLUSION AND SUGGESTIONS FOR FURTHER WORK

The next article will continue to study the asymptotic behavior of some arithmetic functions.

5. ACKNOWLEDGEMENTS

Thanks to everyone who has contributed to the discussion of this paper. I am grateful to everyone who expressed their suggestions and comments in the course of this work.